\documentclass[11pt,leqno]{amsart}
\usepackage{amsmath}
\usepackage{amsthm}
\usepackage{amssymb}
\theoremstyle{definition}

\theoremstyle{plain}

\begin{document}
\setcounter{page}{67}
\begin{flushright}
\scriptsize{{\em Journal of the Nigerian Association of Mathematical Physics}
\\
{\em Volume} {\bf 11}, ({\em November} 2007), 67--70 (Original Print)}
\\
\copyright{\em J. of NAMP}
\end{flushright}
\vspace{16mm}

\begin{center}
{\normalsize\bf Quasi-partial sums of the generalized Bernardi integral of certain analytic functions} \\[12mm]
     {\sc K. O. BABALOLA$^{1}$} \\ [8mm]

\begin{minipage}{100mm}
{\small {\sc Abstract.}

$\frac{\rule{3.90in}{0.02in}}{{}}$\\
In this short note we extend a result of Jahangiri and Farahmand \cite{JM} concerning functions of bounded turning to a more general class of functions.
$\frac{{}}{\rule{3.90in}{0.02in}}$
}
\end{minipage}
\end{center}

 \renewcommand{\thefootnote}{}
 \footnotetext{2000 {\it Mathematics Subject Classification.}
            30C45.}
 \footnotetext{$^1$\,Department of Mathematics, University of Ilorin, Ilorin, Nigeria. ummusalamah.kob\symbol{64}unilorin.edu.ng}

\def\iff{if and only if }
\def\S{Smarandache }
\newcommand{\norm}[1]{\left\Vert#1\right\Vert}

\vskip 12mm

{\bf 1. Introduction }
\medskip

Let $C$ be the complex plane. Denote by $A$ the class of functions:
$$f(z)=z+a_2z^2+\cdots\eqno{(1.1)}$$
which are analytic in the unit disk $E=\{z\colon|z|<1\}$.

In \cite{JM} Jahangiri and Farahmand studied the partial sums of the Liberal integral of the class $B(\beta)$, which consists of functions in $A$ satisfying Re $f'(z)>\beta$, $0\leq\beta<1$. Functions in $B(\beta)$ are called functions of bounded turning. It is known that functions of bounded turning are generally univalent and close-to-convex in the unit disk. In particular they proved that the $mth$ partial sums
$$F_m(z)=z+\sum_{k=2}^m\frac{2}{k+1}a_kz^k\eqno{(1.2)}$$
of the Libera integral
$$F(z)=\int_0^zf(t)dt\eqno{(1.3)}$$
is also of bounded turning. Their result was stated as:\vskip 2mm

{\bf Theorem A}\vskip 2mm

{\em If $\frac{1}{4}\leq\beta<1$ and $f\in B(\beta)$, then $F_m\in B(\frac{4\beta-1}{3})$.}\vskip 2mm

Earlier, Li and Owa \cite{JL} have proved that if $f\in A$ is univalent in $E$, then the partial sums $F_m(z)$ is starlike in the subdisk $|z|<\frac{3}{8}$, the number $\frac{3}{8}$ being the best possible.

The result of Jahangiri and Farahmand \cite{JM} naturally leads to inquistion about a more general class of functions (including $B(\beta)$ as a special case), which was introduced in \cite{TO} by Opoola, and has been studied extensively in \cite{KB}. This is the class $T_n^\alpha(\beta)$ consisting of functions $f\in A$ which satisfy the inequality:
$$Re\frac{D^nf(z)^\alpha}{\alpha^nz^\alpha}>\beta\eqno{(1.4)}$$
where $\alpha>0$ is real, $0\leq\beta<1$, $D^n(n\in N_0=\{0,1,2,\cdots\})$ is the Salagean derivative operator defined as:
$$D^nf(z)=D[D^{n-1}f(z)]=z[D^{n-1}f(z)]'\eqno{(1.5)}$$
with $D^0f(z)=f(z)$ and powers in (1.4) meaning principal values only. Observe that the geometric condition (1.4) slightly modifies the one given originally in \cite{TO} (see \cite{KB}). Observe also that the class $B(\beta)$ corresponds to $n=\alpha=1$.\vskip 2mm

In a recent work we considered the generalized Bernardi integral operator given by: 
$$F(z)^\alpha=\frac{\alpha+c}{z^c}\int_0^zt^{c-1}f(t)^\alpha dt,\;\;\alpha+c>0\eqno{(1.6)}$$ 
and sharpened and extended many earlier results concerning closure, under the integral, of several classes of functions. In the present paper we define a concept of quasi-partial sums and follow a method of Jahangiri and Farahmand \cite{JM} to extend their result (Theorem A) to the class $T_n^\alpha(\beta)$.\vskip 2mm

As we noted in \cite{KO}, the binomial expansion of (1.1) gives
$$f(z)^\alpha=z^\alpha+\sum_{k=2}^\infty a_k(\alpha)z^{\alpha+k-1}\eqno{(1.7)}$$
where $a_k(\alpha)$ is a polynomial depending on the coefficients of $f(z)$ and the index $\alpha$. Hence
$$F(z)^\alpha=z^\alpha+\sum_{k=2}^\infty \frac{\alpha+c}{\alpha+k+c-1}a_k(\alpha)z^{\alpha+k-1}\eqno{(1.8)}$$
and we define the $mth$ quasi-partial sums of the integral (1.6) as follows
$$F_m(z)^\alpha=z^\alpha+\sum_{k=2}^m \frac{\alpha+c}{\alpha+k+c-1}a_k(\alpha)z^{\alpha+k-1}\eqno{(1.9)}$$

In the next section we state the preliminary results.

 \medskip

{\bf 2.0 Preliminary Results}\vskip 2mm

We will require the following lemmas.\vskip 2mm

{\bf Lemma 2.1}(\cite{GG})\vskip 2mm

{\em Let $\theta$ be a real number and $l$ a positive integer. If $-1<\gamma\leq A$, then $\frac{1}{1+\gamma}+\sum_{k=1}^l\frac{\cos k\theta}{k+\gamma}\geq 0$. The constant $A=4.5678018\cdots$ is the best possible.}\vskip 2mm

{\bf Lemma 2.2}\vskip 2mm

{\em For $z\in E$ and $-1<\gamma\leq A=4.5678018\cdots$, $Re\left(\sum_{k=1}^l\frac{z^k}{k+\gamma}\right)\geq-\frac{1}{1+\gamma}$.}

\begin{proof}
Let $z=re^{i\theta}$ where $0\leq r<1$ and $0<|\theta|\leq\pi$. Then by De Moivre's law and the minimum principle for harmonic functions $$Re\left(\sum_{k=1}^l\frac{z^k}{k+\gamma}\right)=\sum_{k=1}^l\frac{r^k\cos k\theta}{k+\gamma}>\sum_{k=1}^l\frac{\cos k\theta}{k+\gamma}.$$
Hence by Abel's lemma \cite{EC} and Lemma 2.1 above the conclusion follows.
\end{proof}

Let $P$ denote the class of functions of the form 
$$p(z)=1+c_1z+\cdots\eqno{(2.1)}$$
normalized by $p(0)=1$ and satisfy Re $p(z)>0$ in $E$. The next lemma concerns convolution of analytic functions with functions in $P$. The convolution (or Hadamard product) of two power series $f(z)=\sum_{k=0}^\infty a_kz^k$ and $g(z)=\sum_{k=0}^\infty b_kz^k$ (written as $f\ast g$) is defined as $(f\ast g)(z)=\sum_{k=0}^\infty a_kb_kz^k$.\vskip 2mm

{\bf Lemma 2.3}(\cite{AW})\vskip 2mm

{\em Let $p(z)$ be analytic in $E$ and satisfy $p(0)=1$ and Re $p(z)>\frac{1}{2}$ in $E$. For analytic function $q(z)$ in $E$, the convolution $p\ast q$ takes values in the convex hull of the image of $E$ under $q(z)$.}

\medskip

{\bf 3.0 Main Results}\vskip 2mm

{\bf Theorem 3.1}\vskip 2mm

{\em Let $f(z)$ given by (1.1) be in the class $T_n^\alpha(\beta)$. Then $$Re\frac{D^nF_m(z)^\alpha}{\alpha^nz^\alpha}>1-\frac{2(1-\beta)(\alpha+c)}{\alpha+c+1},\;\;\alpha+c\leq4.5678018\cdots\eqno{(3.1)}$$. Furthermore if $\beta\geq\frac{1}{2}\frac{\alpha+c-1}{\alpha+c}$, then $F_m(z)$ belongs to some subclasses of the class $T_n^\alpha(\beta)$.}

\begin{proof}
From (1.7) and the condition (1.4) we have
$$Re\left\{1+\frac{1}{2(1-\beta)}\sum_{k=2}^\infty\left(\frac{\alpha+k-1}{\alpha}\right)^na_k(\alpha)z^{k-1}\right\}>\frac{1}{2}\eqno{(3.2)}$$
Also from (1.9) we have
$$\frac{D^nF_m(z)^\alpha}{\alpha^nz^\alpha}=1+\sum_{k=2}^\infty\left(\frac{\alpha+k-1}{\alpha}\right)^n\frac{\alpha+c}{\alpha+c+k-1}a_k(\alpha)z^{k-1}=p\ast q\eqno{(3.3)}$$
where
$$p(z)=1+\frac{1}{2(1-\beta)}\sum_{k=2}^\infty\left(\frac{\alpha+k-1}{\alpha}\right)^na_k(\alpha)z^{k-1},\eqno{(3.4)}$$
$$q(z)=1+2(1-\beta)\sum_{k=2}^\infty\frac{\alpha+c}{\alpha+c+k-1}z^{k-1}.\eqno{(3.5)}$$
Thus by Lemma 2.3 and the condition (3.1), the geometric quantity $\frac{D^nF_m(z)^\alpha}{\alpha^nz^\alpha}$ takes values in the convex hull of $q(E)$. But Re $$Re\;q(z)=1+2(1-\beta)(\alpha+c)Re\left(\sum_{k=1}^\infty\frac{z^k}{\alpha+c+k}\right)\eqno{(3.6)}.$$ We know from (1.6) that $\alpha+c>0$. Now suppose $\alpha+c\leq4.5678018\cdots$, then by taking $l=m-1$ in Lemma 2.2, the real part of the series on the right of (3.6) is greater than $-(\alpha+c+1)^{-1}$ so that 
$$Re\frac{D^nF_m(z)^\alpha}{\alpha^nz^\alpha}=Re q(z)>1-\frac{2(1-\beta)(\alpha+c)}{\alpha+c+1}.\eqno{(3.7)}$$

Now observe that the real number $1-\frac{2(1-\beta)(\alpha+c)}{\alpha+c+1}$ is nonnegative only for $\beta\geq\frac{1}{2}\frac{\alpha+c-1}{\alpha+c}$. Thus only in this case it is clear $F_m(z)$ belongs to some subclasses of the class $T_n^\alpha(\beta)$. This completes the proof.
\end{proof}

\medskip

{\bf Remark} For $\alpha=1$, $c=0$, the partial sums $$F_m(z)=z+\sum_{k=2}^m\frac{a_k}{k}z^k\eqno{(3.8)}$$ of the integral $$F(z)=\int_0^zt^{-1}f(t)dt\eqno{(3.9)}$$ for each $f\in B_n(1)$, belongs to the class $B_n(1)$ in general. More particularly, the partial sum (3.8) of the integral (3.9) of a function of bounded turning in the unit disk is also a function of bounded turning in the unit disk.
 
\medskip

{\bf 4.0 Conclusion}\vskip 2mm

In the paper we defined a new concept of quasi-partial sums of the generalized Bernardi integral. We used the new concept to extend an earlier result of Jahangiri and Farahmand \cite{JM} concerning functions of bounded turning to a more general class of functions.

\end{document}